\documentclass[12pt,letterpaper]{article}
\usepackage{t1enc}
\usepackage[latin1]{inputenc}
\usepackage[english]{babel}
\usepackage{graphics}
\begin{document}

\title{ Sharp inequalities on circular and hyperbolic functions using a Bernoulli inequality type}
\author{Abd Raouf Chouikha
\footnote
{chouikha@math.univ-paris13.fr. 4, Cour des Quesblais 35430 Saint-Pere   
}}
\date{}
\maketitle
\bigskip

\begin{abstract}

In this paper, new sharp bounds for circular and hyperbolic functions are proved. We provide some improvements of previous results by using infinite products, power series expansions and a variant of the so-called Bernoulli inequality. New proofs, refinements as well as new results are offered.


\end{abstract}

\bigskip

During the past several years, sharp inequalities involving circular and hyperbolic functions
have received a lot of attention. Thanks to their usefulness in all areas of mathematics. Old and new
such inequalities, as well as refinements of the so-called Jordan's, Cusa-Huygens and Wilker inequalities,
can be found in [5], and the references therein.\\
In this paper we present a new variant of the standard Bernoulli inequality. This is more general than the one proved in [1]. It permits us to deduce some bounds for  circular and hyperbolic functions. These bounds appear to be sharper than those proved previously.\\ 
In this paper we provide new lower bounds for the function $\cos(x)$ and $\frac{\sin(x)}{x}$  as well as bounds for the products $\frac{\sinh(x) \sin(x)}{x^2}$ and $\cosh(x) \cos(x),$ improving some of those established in the literature.

\section{ A variant of the Bernoulli inequality}
The following two results concerning circular functions, they have been proved by C. Chesneau and Y. Bagul [1]. They improve theorems 1,2 of [2].\\
The first result concerns bounds for $\cos(x)$.

\bigskip
{\bf Proposition 2}\quad {\it For $\alpha \in (0,\pi/2)$ and $x \in (0,\alpha)$, we have 
$$e^{-\beta x^2} \leq \cos (x) \leq e^{-x^2/2},$$
with $\beta = -\frac{log(cos(x\alpha))}{\alpha^2}.$}
\bigskip

The second concerns bounds for $\frac {\sin(x)}{x}$.

\bigskip
{\bf Proposition 3}\quad {\it For $\alpha \in (0,\pi/2)$ and $x \in (0,\alpha)$, we have 
$$e^{-\gamma x^2} \leq \cos (x) \leq e^{-x^2/2},$$
with $\gamma = -\frac{log(sin(x\alpha))}{\alpha^2}.$}
\bigskip

Their proofs are based on infinite products of circular functions as well as the following standard Bernoulli inequality

\bigskip
{\bf Proposition 1}\quad {\it For $u,v \in (0,1)$, we have 
$$1 - uv \geq (1 - v)^u.$$}
\bigskip

In this paper we propose at first the following main result which improves the preceding

\bigskip
{\bf Theorem 1-1}\quad {\it For $u,v \in (0,1)$, we have 
$$1 - uv \geq (1 - v)^{u^2} e^{uv(u-1)} \geq (1 - v)^u.$$}
\bigskip

{\it Proof of Theorem 1-1} \quad For \ $u,v \in(0,1) $ and $k \geq 1$,  one has $u^k \leq u$

 The logarithm series expansion gives 
$$\log (1-uv) = - \sum_{k\geq 1} \frac{u^k v^k}{k} = - uv -\sum_{k\geq 2} \frac{u^k v^k}{k}$$ $$ \geq - uv + u^2 (- \sum_{k\geq 2} \frac{v^k}{k}) = -uv + u^2 (v + \log (1-v).$$
Composing by the exponential function, we obtain the left inequality. To prove the right inequality it suffices to note that the function
$$f(v) = (1 - v) e^v $$ is bounded by $1$ for $v \in (0,1)$. Then powering by $u - u^2 \geq 0$ one gets 
$$ f(v)^{u-u^2} = (1-v)^{u-u^2} e^{v(u-u^2)} \leq 1$$
Implying $$(1-v)^{u^2-u} e^{v(u^2-u)} \geq 1 $$ or equivalently $$(1 - v)^{u^2} e^{uv(u-1)} \geq (1 - v)^u$$ since $u^2 \leq u$.\\

Theorem 1-1 is a sort of a refinement of the classical Bernoulli inequality. It allows us in particular to improve the lower bounds of preceding Propositions 1 and 2. It also allows to bring many other inequalities as well as with sharp bounds.

\bigskip
{\bf Proposition 1-2}\quad {\it For $\alpha \in (0,\pi/2)$ and $x \in (0,\alpha)$, we have the following inequalities}
$$(\cos (\alpha))^{x^2} \leq (\cos (\alpha))^{x^4} e^{(\frac{x^4-\alpha^2 x^2}{2\alpha^2})} \leq \cos(x) \leq e^{-x^2/2}$$
\bigskip

{\it Proof of Proposition 1-2} \quad We will use the infinite product of the cosine function. For $x \in IR$ we have
$$\cos (x) = \prod_{k\geq 1} ( 1 - \frac{4 x^2}{\pi^2(2k-1)^2})$$
The upper bound has been proved (Proposition 2 of [1]). The lower bound follows from Theorem 1. We have to prove the middle bound. With respect that $x \in (0,\alpha)$ using the infinite product and Theorem 1 we may write following [1] 
$$\cos (x) = \prod_{k\geq 1} ( 1 - \frac{4 \alpha^2 x^2}{\alpha^2 \pi^2(2k-1)^2}) = \prod_{k\geq 1} ( 1 - \frac{4 \alpha^2 }{\pi^2(2k-1)^2} \frac{x^2}{\alpha^2})$$ 
$${\Large \geq \prod_{k\geq 1} ( 1 - \frac{4 \alpha^2 x^2}{\pi^2(2k-1)^2})^{\frac{x^4}{\alpha^4}} \ e^{\frac{4 \alpha^2 }{\pi^2(2k-1)^2} (\frac{x^4}{\alpha^4} - \frac{x^2}{\alpha^2})}}$$
$$= {\Large (\cos(\alpha))^{\frac{x^4}{\alpha^4}}\ e^{\sum_{k\geq 1} (\frac{4 \alpha^2 }{\pi^2(2k-1)^2} (\frac{x^4}{\alpha^4} - \frac{x^2}{\alpha^2}))}}.$$

On the other hand the sum may be written $$\sum_{k\geq 1} \frac{4 \alpha^2 }{\pi^2(2k-1)^2} (\frac{x^4}{\alpha^4} - \frac{x^2}{\alpha^2}) = (\frac{x^4}{\alpha^4} - \frac{x^2}{\alpha^2}) \frac{4 \alpha^2 }{\pi^2} \sum_{k\geq 1} \frac{1}{(2k-1)^2}$$ $$= (\frac{x^4}{\alpha^4} - \frac{x^2}{\alpha^2}) \frac{4 \alpha^2 }{\pi^2}{\pi^2} \frac{\pi^2}{8} = (\frac{x^4}{\alpha^4} - \frac{x^2}{\alpha^2}) \frac{ \alpha^2}{2}$$
since $\sum_{k\geq 1} \frac{1}{(2k-1)^2} = \frac{\pi^2}{8}.$\\

Finally, we obtain the inequality 
$$ \cos (x) \geq (\cos(\alpha))^{x^4} \ e^{(\frac{x^4}{\alpha^4} - \frac{x^2}{\alpha^2}) \frac{ \alpha^2 }{2}}.$$\\

The result below gives a generalization of [1, Proposition 3 ]. It improves the lower bound\\

\bigskip
{\bf Proposition 1-3}\quad {\it For $\alpha \in (0,\pi/2)$ and $x \in (0,\alpha)$, we have the following inequalities }
$$(\frac{\sin (\alpha)}{\alpha})^{\frac{x^2}{\alpha^2}} \leq (\frac{\sin (\alpha)}{\alpha})^{\frac{x^4}{\alpha^4}}\ e^{(\frac{x^4-\alpha^2 x^2}{6\alpha^2})} \leq \frac{\sin(x)}{x} \leq e^{-x^2/6}$$
\bigskip

{\it Proof of Proposition 1-3} \quad We will use the infinite product of the sine function the so-called Euler expansion. For $x \in IR$ we have
$$\frac{\sin (x)}{x} = \prod_{k\geq 1} ( 1 - \frac{ x^2}{\pi^2 k^2})$$
The upper bound has been proved (Proposition 2 of [1]). The lower bound follows from Theorem 1. We have to prove the middle bound. With respect that $x \in (0,\alpha)$ using the infinite product and Theorem 1 we may write following [1] 
$$\frac{\sin (x)}{x} = \prod_{k\geq 1} ( 1 - \frac{ \alpha^2 x^2}{\alpha^2 \pi^2 k^2}) = \prod_{k\geq 1} ( 1 - \frac{ \alpha^2 }{\pi^2 k^2} \frac{x^2}{\alpha^2})$$ 
$${\Large \geq \prod_{k\geq 1} ( 1 - \frac{ \alpha^2 x^2}{\pi^2 k^2})^{\frac{x^4}{\alpha^4}} \ e^{\frac{ \alpha^2 }{\pi^2 k^2} (\frac{x^4}{\alpha^4} - \frac{x^2}{\alpha^2})}}$$
$$= {\Large (\frac{\sin(\alpha)}{\alpha})^{\frac{x^4}{\alpha^4}}\ e^{\sum_{k\geq 1} (\frac{ \alpha^2}{\pi^2 k^2} (\frac{x^4}{\alpha^4} - \frac{x^2}{\alpha^2}))}}.$$

On the other hand the sum may be written $$\sum_{k\geq 1} \frac{ \alpha^2 }{\pi^2 k^2} (\frac{x^4}{\alpha^4} - \frac{x^2}{\alpha^2}) = (\frac{x^4}{\alpha^4} - \frac{x^2}{\alpha^2}) \frac{ \alpha^2 }{\pi^2} \sum_{k\geq 1} \frac{1}{ k^2}$$ $$= (\frac{x^4}{\alpha^4} - \frac{x^2}{\alpha^2}) \frac{ \alpha^2 }{\pi^2} \frac{\pi^2}{6} = (\frac{x^4}{\alpha^4} - \frac{x^2}{\alpha^2}) \frac{ \alpha^2 }{6}$$
since $\sum_{k\geq 1} \frac{1}{k^2} = \frac{\pi^2}{6}.$\\

Finally, we obtain the inequality 
$$ \frac{\sin (x)}{x} \geq (\frac{\sin(\alpha)}{\alpha})^{\frac{x^4}{\alpha^4}} \ e^{(\frac{x^4}{\alpha^4} - \frac{x^2}{\alpha^2}) \frac{ \alpha^2}{6}}.$$\\

\section{Others applications on the lower bounds}
In [3, Theorem 1.26] the authors proved for $x\in (0,\pi)$ the following inequality 
$$\frac{sin(x)}{x} \geq (1 - \frac{x^2}{\pi^2})^{\frac{\pi^2}{6}}.$$
The following improves this lower bound

\bigskip
{\bf Proposition 2-4}\quad {\it For $x \in (0,\pi)$ we have the inequalities}
$$\frac{sin(x)}{x} \geq (1 - \frac{x^2}{\pi^2})^{\frac{\pi^4}{90}}\ e^{{x^2}(\frac{\pi^2}{90}-\frac{1}{6})} \geq (1 - \frac{x^2}{\pi^2})^{\frac{\pi^2}{6}}.$$
\bigskip

{\it Proof of Proposition 2-4} \ By Theorem 1 and the infinite product we get
$$\frac{\sin (x)}{x} = \prod_{k\geq 1} ( 1 - \frac{ x^2}{\pi^2 k^2}) \geq \prod_{k\geq 1} (1 - \frac{x^2}{\pi^2})^{\frac{1}{k^4}}\ e^{\frac{x^2}{\pi^2 k^2}(\frac{1}{k^2}-1)}$$
$$= (1 - \frac{x^2}{\pi^2})^{\sum_{k\geq 1}}\frac{1}{k^4}\ e^{\frac{x^2}{\pi^2}(\sum_{k\geq 1}\frac{1}{k^4}-\sum_{k\geq 1}\frac{1}{k^2})}$$
$$ = (1 - \frac{x^2}{\pi^2})^{\frac{\pi^4}{90}}\ e^{\frac{x^2}{\pi^2}(\frac{\pi^4}{90}-\frac{\pi^2}{6})}.$$\\

On the other hand, the well known inequality of the product \\ $\sin(x) \sinh(x) \leq x^2$ \ has been improved in [4, Propositions 2.1 and 2.2] in using infinite product and Bernoulli inequality. We propose the following which is a refinement than the previous.

\bigskip
{\bf Proposition 2-5}\quad {\it For $x \in (0,\pi)$ we have the inequalities}
$$ (\frac{\sin(\alpha) \sinh(\alpha)}{\alpha^2})^\frac{x^4}{\alpha^4} \leq  (\frac{\sin(\alpha) \sinh(\alpha)}{\alpha^2})^\frac{x^8}{\alpha^8}\ e^{[\frac{x^4}{90 }(\frac{x^4}{\alpha^4} - 1)]} 
\leq \frac{\sin(x) \sinh(x)}{x^2} \leq  e^{-\frac{x^4}{90}}$$
$$x^2  (1 - \frac{x^4}{\pi^4})^{\frac{\pi^4}{90}} \leq x^2  (1 - \frac{x^4}{\pi^4})^{\frac{\pi^8}{9450}} \ e^{[\frac{x^4}{90 }(\frac{x^4}{\pi^4}-1)]} \leq  \frac{\sin(x) \sinh(x)}{x^2}$$
\bigskip                                                                                                            
        
{\it Proof of Proposition 2-5} \ The upper bound has been proved (Proposition 2.1 of [4]). The lower bound follows from Theorem 1. We have to prove the middle bound. With respect that $x \in (0,\alpha)$ using the infinite products
$$\frac{\sin (x)}{x} = \prod_{k\geq 1} ( 1 - \frac{ x^2}{\pi^2 k^2}), \qquad \frac{\sinh (x)}{x} = \prod_{k\geq 1} ( 1 + \frac{ x^2}{\pi^2 k^2}).$$
Therefore $$\sin(x) \sinh(x) = x^2  \prod_{k\geq 1} ( 1 - \frac{ x^4}{\pi^4 k^4}).$$
For this bound we will use Theorem 1. It follows 
$$\sin(x) \sinh(x) = x^2  \prod_{k\geq 1} ( 1 - \frac{\alpha^4 x^4}{\pi^4 \alpha^4 k^4})$$
$$\geq \prod_{k\geq 1} (1 - \frac{\alpha^4}{\pi^4 k^4})^\frac{x^8}{\alpha^8}\ e^{[\frac{x^4}{\pi^4 k^4}(\frac{x^4}{\alpha^4} - 1)]}$$
$$ = x^2 (\frac{\sin(\alpha) \sinh(\alpha)}{\alpha^2})^\frac{x^8}{\alpha^8} \ e^{[{\frac{x^4}{\pi^4 }(\frac{x^4}{\alpha^4} - 1)\sum_{k\geq 1}\frac{1}{k^4}}]}.$$
It implies the inequality since $\sum_{k\geq 1} \frac{1}{k^4} = \frac{\pi^4}{90}.$\\
For the second inequality which is sharp than the preceding we use again the infinite products
$$\sin(x) \sinh(x) = x^2  \prod_{k\geq 1} ( 1 - \frac{ x^4}{\pi^4 k^4}).$$
By Theorem 1 we may write 
 $$\sin(x) \sinh(x) =x^2  \prod_{k\geq 1} ( 1 - \frac{ x^4}{\pi^4}\frac{1}{k^4})  \geq x^2  \prod_{k\geq 1}(1 - \frac{x^4}{\pi^4})^{\frac{1}{k^8}} \ e^{[\frac{x^4}{\pi^4 k^4}(\frac{x^4}{\pi^4}-1)]}$$
$$= x^2  (1 - \frac{x^4}{\pi^4})^{\zeta(8)} \ e^{\frac{x^4}{\pi^4 }(\frac{x^4}{\pi^4}-1)\zeta(4)}.$$
Thus $$\sin(x) \sinh(x) \geq x^2  (1 - \frac{x^4}{\pi^4})^{\frac{\pi^8}{9450}} \ e^{[\frac{x^4}{\pi^4 }(\frac{x^4}{\pi^4}-1)\frac{\pi^4}{90}]}.$$\\
{
Similarly, for the cosine product $\cos(x) \cosh(x)$ we have the following result which determines a double inequality improving [4, Propositions 2.3 and 2.4]

\bigskip
{\bf Proposition 2-6}\quad {\it For $x \in (0,\pi)$ we have the inequalities}
$$ [\cos(\alpha) \cosh(\alpha)]^\frac{x^4}{\alpha^4} \leq  [\cos(\alpha) \cosh(\alpha)]^\frac{x^8}{\alpha^8}\ e^{[{\frac{16 x^4}{96 }(\frac{x^4}{\alpha^4} - 1)}]} 
\leq \cos(x) \cosh(x) \leq  e^{-\frac{x^4}{6}}$$
$$(1 - \frac{16 x^4}{\pi^4})^{\frac{\pi^4}{90}} \leq  (1 - \frac{16 x^4}{\pi^4})^{\frac{17\pi^8}{161280}} \ e^{[\frac{16 x^4}{96 }(\frac{x^4}{\pi^4}-1)]} \leq \cos(x) \cosh(x).$$
\bigskip

 {\it Proof of Proposition 2-6} \ The upper bound has been proved (Proposition 2.3 of [4]). The lower bound follows from Theorem 1. We have to prove the middle bound. With respect that $x \in (0,\alpha)$ using the infinite products
For $x \in IR$ we have
$$\cos (x) = \prod_{k\geq 1} ( 1 - \frac{4 x^2}{\pi^2 (2k-1)^2}), \qquad \frac{\cosh (x)}{x} = \prod_{k\geq 1} ( 1 + \frac{4 x^2}{\pi^2 (2k-1)^2}).$$
Therefore $$\cos(x) \cosh(x) =   \prod_{k\geq 1} ( 1 - \frac{16 x^4}{\pi^4 (2k-1)^4}).$$
For this bound we will use Theorem 1. It follows 
$$\cos(x) \cosh(x) =   \prod_{k\geq 1} ( 1 - \frac{16\alpha^4 x^4}{\pi^4 \alpha^4 (2k-1)^4})$$
$$\geq \prod_{k\geq 1} [1 - \frac{\alpha^4}{\pi^4 (2k-1)^4}]^\frac{x^8}{\alpha^8}\ e^{[\frac{16 x^4}{\pi^4 (2k-1)^4}(\frac{x^4}{\alpha^4} - 1)]}$$
$$ =  [\cos(\alpha) \cosh(\alpha)]^\frac{x^8}{\alpha^8} \ e^{[{\frac{16 x^4}{\pi^4 }(\frac{x^4}{\alpha^4} - 1)\sum_{k\geq 1}\frac{1}{(2k-1)^4}}]}.$$
It implies the inequality since $\sum_{k\geq 1} \frac{1}{(2k-1)^4} = \frac{\pi^4}{96}.$\\
For the second inequality which is sharp than the preceding we use again the infinite products
$$\cos(x) \cosh(x) =   \prod_{k\geq 1} ( 1 - \frac{16 x^4}{\pi^4 (2k-1)^4}).$$
By Theorem 1 we may write 
 $$\cos(x) \cosh(x) =  \prod_{k\geq 1} ( 1 - \frac{16 x^4}{\pi^4}\frac{1}{(2k-1)^4})  \geq x^2  \prod_{k\geq 1}[1 - \frac{16 x^4}{\pi^4}]^{\frac{1}{(2k-1)^8}} \ e^{[\frac{16 x^4}{\pi^4 (2k-1)^4}(\frac{x^4}{\pi^4}-1)]}$$
$$=   [1 - \frac{16 x^4}{\pi^4}]^{\frac{17\pi^8}{161280}} \ e^{[\frac{16 x^4}{\pi^4 }(\frac{x^4}{\pi^4}-1)\frac{\pi^4}{96}]}$$
since $ \sum_{k\geq 1} \frac{1}{(2k-1)^8} =  \frac{17\pi^8}{161280}$.\\
Thus $$\cos(x) \cosh(x) \geq [1 - \frac{16 x^4}{\pi^4}]^{\frac{17\pi^8}{161280}} \ e^{[\frac{x^4}{96 }(\frac{x^4}{\pi^4}-1)]}.$$

\vspace{3cm}

{\bf REFERENCES }\\

[1] \ C. Chesneau and Y. Bagul, \ {\it A note on some new bounds for trigonometric functions using infinite products}, 2018. hal-01934571v2.\\

[2] \ Y. Bagul, \ {\it Inequalities involving circular, hyperbolic and exponential functions}, J. of Math. Ineq, vol 11, n 3, pp 695-699, 2017.\\

[3] \ B. Bhayo and J. Sandor, \ {\it On Jordan's, Redheffer's and Wilker's inequality}, J. of Math. Ineq, vol 19, n 3, pp 823-839, 2016.\\ 

[4] \ C. Chesneau and Y. Bagul, \ {\it New refinements of two well-known inequalities}, https://www.researchgate.net/publication/332849515.\\

[5] \ L. Zhu, \ {\it A source of inequalities for circular functions}, Comp. and Math. with Appl., 58, pp 1998-2004, 2009.\\

\end{document}